\documentclass[12pt]{article}
\usepackage{amsmath,amsfonts,amssymb,rotating,colordvi}
\usepackage{ytableau}
\usepackage{color}
\usepackage[unicode=true,
 bookmarks=true,bookmarksnumbered=true,bookmarksopen=true,bookmarksopenlevel=2,
 breaklinks=false,pdfborder={0 0 1},backref=false,colorlinks=true]
 {hyperref}

\hypersetup{pdftitle={},
 pdfauthor={},
 pdfsubject={ },
 pdfkeywords={ },
 linkcolor=black,citecolor=black,urlcolor=blue,filecolor=blue,
pdfpagelayout=OneColumn,pdfnewwindow=true,pdfstartview=XYZ,plainpages=false}

\def\qed{\nopagebreak\hfill{\rule{4pt}{7pt}}}
\def\proof{\noindent {\it{Proof.} \hskip 2pt}}

 \newtheorem{thm}{Theorem}[section]

\newtheorem{lem}[thm]{Lemma}

\newtheorem{coro}[thm]{Corollary}
\newtheorem{conj}[thm]{Conjecture}

\numberwithin{equation}{section}
\DeclareMathOperator{\asc}{asc}

\DeclareMathOperator{\des}{des}

\newdimen\Squaresize \Squaresize=11pt
\newdimen\Thickness \Thickness=0.7pt
\def\Square#1{\hbox{\vrule width \Thickness
   \vbox to \Squaresize{\hrule height \Thickness\vss
    \hbox to \Squaresize{\hss#1\hss}
   \vss\hrule height\Thickness}
\unskip\vrule width \Thickness} \kern-\Thickness}

\def\Vsquare#1{\vbox{\Square{$#1$}}\kern-\Thickness}

\def\moins{\raise 1pt\hbox{{$\scriptstyle -$}}}

\begin{document}

\begin{center}
\textbf{\large{}Kirillov's unimodality conjecture for the rectangular Narayana
polynomials}\textbf{ }
\par\end{center}

\begin{center}
Herman Z.Q. Chen$^{1}$, Arthur L.B. Yang$^{2}$, Philip B. Zhang$^{3}$\\[6pt]
\par\end{center}

\begin{center}
$^{1,2}$Center for Combinatorics, LPMC\\
 Nankai University, Tianjin 300071, P. R. China

$^{3}$College of Mathematical Science \\
Tianjin Normal University, Tianjin  300387, P. R. China\\[6pt]
\par\end{center}

\begin{center}
Email: $^{1}$\texttt{zqchern@163.com, $^{2}$\texttt{yang@nankai.edu.cn}},
$^{3}$\texttt{zhangbiaonk@163.com}
\par
\end{center}

\noindent \emph{Abstract.}
In the study of Kostka numbers and Catalan
numbers, Kirillov posed a unimodality conjecture for the rectangular Narayana
polynomials.
We prove that the rectangular Narayana polynomials have only real zeros, and
thereby confirm Kirillov's unimodality conjecture with the help of Newton's
inequality.
By using an equidistribution property between descent numbers and ascent
numbers
on ballot paths due to Sulanke and a bijection between lattice words and
standard Young tableaux, we show that the rectangular Narayana polynomial is
equal to the descent generating function on standard Young tableaux of certain
rectangular shape, up to a power of the indeterminate.
Then we obtain the real-rootedness of the rectangular Narayana polynomial based
on
Brenti's result that the descent generating function of standard Young tableaux
has only real zeros.

\noindent \emph{AMS Classification 2010:} 05A20, 30C15


\section{Introduction}

The main objective of this paper is to prove a unimodality conjecture for the
rectangular Narayana polynomials
in the study of Kostka numbers and Catalan numbers. This conjecture was first
posed by Kirillov \cite{Kirillov99} in 1999, and it was restated by himself
\cite{Kirillov15} in 2015. In this paper we prove that the rectangular
Narayana polynomials have only real zeros, an even stronger result than
Kirillov's conjecture.

Let us begin with an overview of Kirillov's conjecture. Throughout this paper, we abbreviate the vector
$(m,m,\ldots,m)$ with $n$ occurrences of $m$ as $(m^n)$ for any positive
integer $m$ and $n$.
We say that a word $\mathtt{w}=w_1w_2\cdots w_{nm}$ in symbols
$1,2,\ldots, m$ is a {lattice word} of weight $(m^n)$,
if the following
conditions hold:
\begin{itemize}
\item[(a)] each $i$ between $1$ and $m$ occurs exactly $n$ times; and

\item[(b)] for each $1\leq r\leq nm$ and $1\leq i \leq m-1$, the number of
$i$'s in $w_1w_2\cdots w_r$ is not less than the number of $(i+1)$'s.
\end{itemize}

Given a word $\mathtt{w}=w_1w_2\cdots w_p$ of length $p$, we say that $i$ is an
{ascent} of $\mathtt{w}$ if $w_i < w_{i+1}$, and a {descent} of
$\mathtt{w}$ if $w_i >w_{i+1}$. Denote the number of ascents of
$\mathtt{w}$ by  $\asc(\mathtt{w})$, and the number of descents
$\des(\mathtt{w})$. For any $m$
and $n$, the rectangular Narayana polynomial $N(n,m;t)$ is defined by
\begin{align}\label{def-Narayana_Pol}
N(n,m;t)= \sum_{\mathtt{w}\in \mathcal{N}(n,m)} t^{\des(\mathtt{w})},
\end{align}
where $\mathcal{N}(n,m)$ is the set of lattice words
of weight $(m^n)$. Note that $N(n,2;t)$ is the classical Narayana polynomial,
and
$N(n,2;1)$ is the classical Catalan number, see \cite{Kirillov15}. For this
reason,
$N(n,m;1)$ is called the rectangular Catalan number.

Kirillov's conjecture is concerned with the unimodality of the rectangular
Narayana polynomial $N(n,m;t)$.
Recall that a sequence $\{a_0,a_1,\ldots,a_n\}$ of positive real numbers is
said to be {unimodal} if there exists an integer $i\geq0$
such that
$$a_0\leq \cdots \leq a_{i-1}\leq a_i \geq a_{i+1}\geq\cdots\geq a_n,$$
and
it is said to be {log-concave} if, for each $1\leq i\leq n-1$, there holds
$$a_i^2\geq a_{i-1}a_{i+1}.$$
Clearly, for a sequence of positive numbers, its log-concavity implies
unimodality.
Given a polynomial with real coefficients
$$f(t)=\sum_{k=0}^n a_k t^k,$$
we say that it is unimodal (or log-concave) if its coefficient sequence
$\{a_0,a_1,\ldots,a_n\}$ is unimodal (resp. log-concave). Kirillov proposed the
following conjecture.

\begin{conj}\cite[Conjecture 2.5]{Kirillov15}\label{Conj-Kirillov}
For any $m$ and $n$, the rectangular Narayana polynomial $N(n,m;t)$ is unimodal
as a polynomial of $t$.
\end{conj}

In this paper, we give an affirmative answer to the above conjecture.
Instead of directly proving its unimodality, we shall show that the rectangular
Narayana polynomial $N(n,m;t)$ has only real zeros. By the well known Newton's inequality, if a polynomial with nonnegative coefficients has only
real zeros, then its coefficient sequence
must be log-concave and unimodal. Thus, from the real-rootedness of $N(n,m;t)$ we deduce its
log-concavity and unimodality.

The remainder of this paper is organized as follows.
In Section \ref{section-2}, we show that the rectangular Narayana polynomial
$N(n,m;t)$ is equal to the descent generating function on standard Young
tableaux of shape $(n^m)$, up to a power of $t$. We use a result of Sulanke
\cite{Sulanke04} that the ascent and descent statistics are equidistributed
over the set of ballot paths.
In Section \ref{section-3}, we first prove the real-rootedness of the descent
generating function on standard Young tableaux, and then obtain the
real-rootedness of $N(n,m;t)$. The key to this approach is a connection between
the descent generating functions of standard Young tableaux and the Eulerian
polynomials of column-strict labeled Ferrers posets. While, the latter
polynomials have only real zeros, as proven by Brenti \cite{Brenti89} in the
study of Neggers-Stanley conjecture.

\section{Tableau interpretation} \label{section-2}

The aim of this section is to interpret the rectangular
Narayana polynomials as the descent generating function on standard
Young tableaux.

Let us first recall some definitions. Given a partition $\lambda$, its
 Young diagram is defined to be an array of squares in the plane
justified from the top left corner with $l(\lambda)$ rows and $\lambda_i$
squares in row i.  By transposing the diagram of $\lambda$, we get the
conjugate partition of $\lambda$, denoted $\lambda'$.
A cell $(i,j)$ of $\lambda$
is  in the $i$-th row from the top and in the $j$-th column from the left. A
{semistandard Young tableau (SSYT)} of shape $\lambda$ is a filling of its
diagram  by positive integers such
that it is weakly increasing in every row and strictly increasing down every
column. The type of $T$ is defined to be
the composition $\alpha=(\alpha_1,\alpha_2,\ldots)$, where $\alpha_i$ is the
number of $i$'s in $T$. If $T$ is of type $\alpha$ with $\alpha_i=1$ for $1\leq
i\leq |\lambda|$ and
$\alpha_i=0$ for $i> |\lambda|$, then it is called a {standard Young tableau
(SYT)} of shape $\lambda$. Let $\mathcal{T}_{\lambda}$ denote the set of SYTs
of shape $\lambda$. Given a standard Young tableau,
we say that $i$ is a {descent} of $T$ if $i+1$ appears in a lower row of $T$
than $i$. Define the descent set $D(T)$ to be the set of all descent of $T$,
and denote by $\des(T)$ the
number of descents of $T$.

The main result of this section is as follows.

\begin{thm} \label{main-thm-0} For any positive integers $m$ and $n$, we have
\begin{align}
N(n,m;t)=t^{1-m}\sum_{T\in \mathcal{T}_{(n^m)}}t^{\des(T)}.
\end{align}
\end{thm}

To prove the above result, we need a bijection between the set of lattice paths
and the set of standard Young tableaux. Here we use a very natural bijection
$\phi$ between the lattice word of weight $(m^n)$
and the
standard Young tableau of shape $(n^m)$,
see \cite[p. 92]{Coleman68}, \cite[p. 221]{Hamermesh89} and \cite{Stanley99}.
To be self-contained, we shall give a description of
this bijection in the following.

Given a lattice word $\mathtt{w}=w_1\cdots w_{nm}$ of weight $(m^n)$,  let
$T=\phi(\mathtt{w})$ be the tableau of shape $(n^m)$ by filling
the  square $(i,j)$ with $k$ provided that $w_k$ is the $j$-th
occurrence of $i$ in $\mathtt{w}$ from left to right. Clearly,
$T$ is a standard Young tableau. Conversely, given a standard Young tableau $T$
of shape $(n^m)$, define a word $\mathtt{w}$ by letting $w_i$ to be
$j$ if $i$ is in the $j$-th row of $T$. It is easy to verify that
$\mathtt{w}=\phi^{-1}(T)$. Figure \ref{fig-1} gives an illustration of this
bijection, where $T$ is of shape $(4^3)$ and  $\mathtt{w}$ is of
weight $(3^4)$.

\begin{figure}[ht]
\begin{center}
$\mathtt{w}=~121113223233$~$\mapsto$~
$T=$~
\begin{ytableau}
 1 & 3 & 4 & 5 \\
 2 & 7 & 8 & 10\\
 6 & 9 & 11 & 12
\end{ytableau}
\end{center}
\caption{Bijection between standard Young tableaux and lattice
words}\label{fig-1}
\end{figure}

By using the above bijection $\phi$, we obtain the following result.

\begin{lem}
For any positive integers $m$ and $n$, we have
\begin{align}\label{D_TandA_w}
 \sum_{T\in \mathcal{T}_{(n^m)}}t^{\des(T)}=\sum_{\mathtt{w}\in\mathcal{N}(n,m)}
t^{\asc(\mathtt{w})}.
\end{align}
\end{lem}

\proof
Suppose that $T=\phi(\mathtt{w})$.
Note that if $i$ is an ascent in $\mathtt{w}$, i.e. $w_i< w_{i+1}$,
then
$i+1$ is filled in the $w_{i+1}$-th row, which is lower than
the row including $i$ in $T$. Thus, $\asc(\mathtt{w})=\des(T)$. This completes
the
proof. \qed

To prove Theorem \ref{main-thm}, it remains to show that
\begin{align}\label{eq-sun}
t^{1-m} \sum_{\mathtt{w}\in\mathcal{N}(n,m)}
t^{\asc(\mathtt{w})}=\sum_{\mathtt{w}\in\mathcal{N}(n,m)}
t^{\des(\mathtt{w})}.
\end{align}
In fact, this has been established by Sulanke \cite{Sulanke04}, which was
stated in terms of ballot paths. In the following, we
shall give an overview of Sulanke's result.

Recall that a ballot path for $m$-candidates is an
$m$-dimensional lattice path running from $(0,0,\ldots,0)$
to $(n,n,\ldots,n)$ with the steps:
\begin{align*}
X_1&:=(1,0,\ldots,0),\\
X_2&:=(0,1,\ldots,0),\\
\vdots &\qquad \qquad  \vdots \\
X_m&:=(0,0,\ldots,1),
\end{align*}
and lying in the region
$$\{ (x_1,x_2,\ldots,x_m):0\leq x_1\leq x_2
\leq\ldots\leq x_m \}.$$ Denote $\mathcal{C}(m,n)$ by the set of all such paths.

For any path $P:=p_1p_2\ldots p_{mn}\in \mathcal{C}(m,n)$,  {the number of
ascents} of $P$ is defined by
$$\asc(P):=|\{i: p_i p_{i+1}=X_j X_l,j<l \}|,$$
and the number of descents of $P$ by
 $$\des(P):=|\{i: p_i p_{i+1}=X_j X_l,j>l \}|.$$

Sulanke \cite{Sulanke04} obtained the following result by a nice bijection.

\begin{lem}\cite[Proposition 2]{Sulanke04}
 For any positive integers $m$ and $n$, we have
 \begin{align}
\sum_{P\in \mathcal{C}(m,n)} t^{\asc(P)}=\sum_{P\in \mathcal{C}(m,n)}
t^{\des(P)-m+1}.
\end{align}
\end{lem}

Note that there is an obvious bijection between $\mathcal{C}(m,n)$ and
$\mathcal{N}(n,m)$: given a path $P\in\mathcal{C}(m,n)$, simply replacing each
step $X_i$ of $P$ by the symbol $m-i+1$, and
the resulting word $\mathtt{w}$ is clearly a lattice word of
$\mathcal{N}(n,m)$. Moreover, we have $\asc(P)=\des(\mathtt{w})$ and
$\des(P)=\asc(\mathtt{w})$.
With this bijection, Sulanke's result can be restated as \eqref{eq-sun}.

\textit{Proof of Theorem \ref{main-thm}.}
Combining \eqref{def-Narayana_Pol}, \eqref{D_TandA_w} and \eqref{eq-sun},
we immediately obtain the desired result. \qed

\section{Real zeros} \label{section-3}

In this section, we aim to prove the real-rootedness of rectangular Narayana
polynomials. Our main result of this section is as follows.

\begin{thm}\label{main-thm}
The rectangular Narayana polynomial $N(n,m;t)$ has only real zeros for any $m$
and $n$.
\end{thm}

By Theorem \ref{main-thm-0}, we only need to show that the following polynomial
$$\sum_{T\in \mathcal{T}_{(n^m)}}t^{\des(T)}$$
has only real zeros. In fact, Brenti \cite{Brenti89} has already obtained a
more general result
during the study of the Neggers-Stanley Conjecture, also known as the Poset
Conjecture.

We now give an overview of the Neggers-Stanley conjecture. Suppose that $P$ is a
finite poset
of cardinality $p$. A labeling $\omega$ of $P$ is a bijection from $P$ to
$\{1,2,\ldots,p\}$. The labeling $\omega$ is called natural if $x\leq y$
implies $\omega(x)\leq \omega(y)$ for any $x,y\in P$.
A {$(P, \omega)$-partition} is a map $\sigma$ which
satisfies the following conditions:
\begin{itemize}
\item[(i)] $\sigma$ is order reversing, namely, $\sigma(x)\geq \sigma(y)$ if
$x\leq y$ in $P$; and moreover
\item[(ii)] if $\omega(x)>\omega(y)$, then $\sigma(x)>\sigma(y)$.
\end{itemize}
The {order polynomial} $\Omega(P,\omega; n)$ is defined as the number of $(P,
\omega)$-partitions
$\sigma$ with $\sigma(x)\leq n$ for any $x\in P$. It is known that
$\Omega(P,\omega; n)$ is a polynomial
of degree $p$ in $n$. By a well known result about rational generating function,
see \cite[Corollary 4.3.1]{Stanley96}, there exists a polynomial
$W(P,\omega;t)$ of degree $\leq p$ such that
\begin{align}\label{Omega_WP}
 \sum_{n\geq 0} \Omega(P,\omega;n+1)t^n=\frac{W(P,\omega;t)}{(1-t)^{p+1}}.
\end{align}
A fundamental result in the theory of $(P, \omega)$-partitions developed by
Stanley \cite{Stanley72} is that the polynomial $W(P,\omega;t)$ can be written
as
\begin{align}
W(P,\omega;t)=\sum_{\pi\in\mathcal{L}(P,\omega)} t^{\des(\pi)},
\end{align}
where $\mathcal{L}(P,\omega)$ is the Jordan-H\"older set
of $(P,\omega)$, that is defined to be the set of
permutations $\pi=\omega(\sigma_1)\omega(\sigma_2)\cdots\omega(\sigma_p)$ where
$\sigma_1\sigma_2\cdots\sigma_p$ is a linear extension of $P$. The polynomials $W(P,\omega;t)$
are  also called $(P,\omega)$\,-Eulerian polynomials. Note that,
when $(P,\omega)$ is a $p$-element anti-chain, the polynomial $W(P,\omega;t)$
is the traditional Eulerian polynomial for the symmetric group on $p$ elements.
The Poset Conjecture is the following.

\begin{conj}\cite[Conjecture 1]{Brenti89} \label{conj-Neggers-Stanely}
For any labeled poset $(P, \omega)$ the polynomial $W(P,\omega;t)$ has only
real zeros as a polynomial of $t$.
\end{conj}

Conjecture \ref{conj-Neggers-Stanely} was formulated for naturally labeled poset
by Neggers \cite{Neggers78} in 1978,
and was generalized to its current form by Stanley in 1986.
It has been proved for some special cases by Brenti \cite{Brenti89},
Wagner \cite{Wagner92}, Reiner \cite{Reiner05} and Br\"{a}nd\'{e}n
\cite{Branden04a}. However,
Br\"{a}nd\'{e}n \cite{Branden04c} and Stembridge \cite{Stembridge07} showed that
the Poset Conjecture doesn't hold in general. One of the most interesting
posets, for which the Poset Conjecture holds, is
the {Ferrers poset} with a column strict labeling.

Recall that the {Ferrers poset $P_{\lambda}$} with respect to
$\lambda=(\lambda_1,\ldots,\lambda_{\ell})$ is the poset
\[ P_{\lambda}=\{(i,j)\in \mathbb{P} \times \mathbb{P}: 1\leq i\leq \ell, 1\leq
j\leq \lambda_i \},\]
ordered by the standard product ordering. We say that a labeling $\omega$ of
$P_{\lambda}$ is {column strict}
if $\omega(i,j)>\omega(i+1,j)$ and $\omega(i,j)<\omega(i,j+1)$ for all
$(i,j)\in P_{\lambda}$.
Note that, given a permutation
$\pi=\pi_1\pi_2\cdots\pi_p\in \mathcal{L}(P_{\lambda},\omega)$, the sequence
$\omega^{-1}(\pi_1)\omega^{-1}(\pi_2)\cdots\omega^{-1}(\pi_p)$ is a linear
extension of $P_{\lambda}$.
Let $T$ be the tableau of shape $\lambda$ by filling
the square $\omega^{-1}(\pi_k)$ with $k$. Clearly,
$T$ is a standard Young tableau. Furthermore,  $k$ is a descent in $\pi$ if and only if
$k$ is a descent in $T$. In fact, suppose that $k$ and $k+1$ are in the square
$(x,y)$ and $(x',y')$ then $\pi_k>\pi_{k+1}$ implies that $x<x'$, that is,
$k+1$ appears in a lower row of $T$ than $k$.
For example, taking $\pi=4215673$ and the labeling $\omega$ showing in the
Figure \ref{perm-tab}, we obtain the standard Young tableau T.

\begin{figure}[ht]
\begin{center}
$\omega=$~
\begin{ytableau}
 4 & 5 & 6 & 7 \\
 2 & 3 \\
 1
\end{ytableau}~,
~~~
$\pi=4215673$~$\mapsto$~
$T=$~
\begin{ytableau}
 1 & 4 & 5 & 6 \\
 2 & 7 \\
 3
\end{ytableau}
\end{center}
\caption{Bijection between permutations in $\mathcal{L}(P_{\lambda},\omega)$ and standard
Young tableaux of shape $\lambda$ for a given labeling $\omega$.}\label{perm-tab}
\end{figure}
Therefore,
\begin{align}
W(P_{\lambda},\omega;t)=\sum_{T\in \mathcal{T}_{\lambda}}t^{\des(T)}.
\end{align}

Brenti \cite{Brenti89} proved the following result, see also Br\"{a}nd\'{e}n
\cite{Branden04}.

\begin{thm}\label{WPrz}\cite[p. 60, Proof of Theorem 5.3.2]{Brenti89}
 Let $(P_{\lambda},\omega)$ be labeled column strict.
 Then $W(P_{\lambda},\omega;t)$ has only real zeros, namely the polynomial
$$\sum_{T\in \mathcal{T}_{\lambda}}t^{\des(T)}$$
has only real zeros.
\end{thm}

Now we can give a proof of Theorem \ref{main-thm}.

\noindent \textit{Proof of Theorem \ref{main-thm}.} This follows
from Theorems
\ref{main-thm-0} and \ref{WPrz}.

As an immediate corollary of Theorem \ref{main-thm}, we obtain the following result, which gives 
an affirmative answer to Kirillov's conjecture.  

\begin{coro}
The rectangular Narayana polynomial $N(n,m;t)$ is unimodal for any $m$ and $n$.
\end{coro}

\noindent{\bf Acknowledgements.} This work was supported by the 973 Project, the PCSIRT Project of the Ministry of Education and the National Science Foundation of China.

\end{document}